\documentclass[12pt]{amsart}
\usepackage[utf8]{inputenc}
\usepackage[T1]{fontenc}
\usepackage{microtype}
\usepackage{fullpage}
\usepackage{cleveref}
\usepackage{enumerate}
\usepackage{amssymb}
\usepackage{amsmath}
\usepackage{amsthm}
\usepackage{mathtools}
\usepackage{mathrsfs}
\usepackage{amsrefs}
\usepackage{tikz}
\usepackage{color}
\newcommand{\R}{\mathbb{R}}
\newcommand{\Z}{\mathbb{Z}}
\renewcommand{\S}{\mathbb{S}}

\newcommand{\cM}{\mathcal{M}}
\newcommand{\cQ}{\mathcal{Q}}
\DeclareMathOperator{\BV}{\dot{BV}}
\newcommand{\loc}{\text{loc}}
\newcommand{\ga}{\gamma}
\newcommand{\la}{\lambda}
\newcommand{\cB}{\mathcal{B}}
\newcommand{\dif}{\mathrm{d}}
\newcommand{\cE}{\mathcal{E}}

\newtheorem{thm}{Theorem}
\title{Sobolev spaces revisited}
\author{Ha\"{i}m Brezis}
\author{Andreas Seeger}
\author{Jean Van Schaftingen}
\author{Po-Lam Yung}

\address{Department of Mathematics\\
  Rutgers University, Hill Center, Busch Campus\\ 
  110 Frelinghuysen Road, Piscataway, NJ 08854, USA}
\address{
  Departments of Mathematics and Computer Science\\ Technion, Israel Institute of Technology\\ 32.000 Haifa, Israel}
\address{
  Laboratoire Jacques-Louis Lions\\
  Sorbonne Universit\'es, UPMC Universit\'e Paris-6, 4  place Jussieu\\
  75005 Paris, France}
\email{brezis@math.rutgers.edu}

\address{Department of Mathematics \\ 
University of Wisconsin, Madison \\
480 Lincoln Drive, Madison, WI, 53706, USA}
\email{seeger@math.wisc.edu}

\address{Universit\'e catholique de Louvain\\ 
Institut de Recherche en Math\'ematique et Physique\\
Chemin du Cyclotron 2 bte L7.01.01\\
1348 Louvain-la-Neuve\\
Belgium}
\email{Jean.VanSchaftingen@UCLouvain.be}

\address{Mathematical Sciences Institute \\
Australian National University \\
Canberra ACT 2601 \\
Australia} 
\email{PoLam.Yung@anu.edu.au}

\dedicatory{Dedicated with emotion to the memory of Antonio Ambrosetti}

\begin{document}

\begin{abstract}
We describe a recent, one-parameter family of characterizations of Sobolev and BV functions on $\R^n$, using sizes of superlevel sets of suitable difference quotients. 
This provides an alternative point of view to the BBM formula by Bourgain, Brezis and Mironescu, and complements in the case of BV  some results of Cohen, Dahmen, Daubechies and DeVore about the sizes of wavelet coefficients of such functions.
An application towards Gagliardo-Nirenberg interpolation inequalities is then given. 
We also establish a related one-parameter family of formulae for the $L^p$ norm of functions in $L^p(\R^n)$.
\end{abstract}

\maketitle

\section{Introduction}

In this note we revisit Sobolev spaces on $\R^n$, $n \geq 1$, from the point of view of difference quotients. For $1 \leq p < \infty$, the homogeneous Sobolev space $\dot{W}^{1,p}$ on $\R^n$ consists of all locally integrable functions $u$ modulo constants, whose distributional gradient $\nabla u \in L^p$. It is normed by
\[
\|\nabla u\|_{L^p} = \Big( \int_{\R^n} |\nabla u|^p \, \dif x \Big)^{1/p}.
\]
We will also work with the homogeneous $\BV$ space (BV stands for \emph{bounded variation}). It consists of all locally integrable functions $u$ modulo constants, whose distributional gradient $\nabla u$ is a finite Radon measure (written $\nabla u \in \cM$). 
In other words, it is the space of all $u \in L^1_{\loc}$ such that
\begin{equation}
\label{eq_Aef7bohJ2pieShooDeupeeli}
\sup \left\{ \Big| \int_{\R^n} u(x) \operatorname{div} \phi(x) \, \dif x \Big| \colon \phi \in C^1_c(\R^n;\R^n), \|\phi\|_{L^{\infty}(\R^n;\R^n)} \leq 1 \right\} 
\end{equation}
is finite (which in particular contains $\dot{W}^{1,1}$). It is normed by 
\[
\|u\|_{\BV} \coloneqq \|\nabla u\|_{\cM}
\] 
which is just the supremum in \eqref{eq_Aef7bohJ2pieShooDeupeeli}.

In \cite{Brezis_VanSchaftingen_Yung_2021}, a new formula was established for $\|\nabla u\|_{L^p(\R^n)}$ for $u \in C^{\infty}_c$ that involves only difference quotients and no gradients.
In \cite{BSVY}, this formula was extended to all $u \in \dot{W}^{1,p}$, or all $u \in \BV$, and in fact we found a natural one-parameter family of such formulae.
A related limiting formula has been obtained earlier by Nguyen \cite{Nguyen06} and by Brezis and Nguyen \cite{BN2018}.
The one-parameter family of formulae we found can be used to recover certain Gagliardo-Nirenberg interpolation inequalities due to Cohen, Dahmen, Daubechies and DeVore \cite{CDDD}.
It also allows us to prove some substitutes when such interpolation inequalities fail (cf.\ \cite{Brezis_VanSchaftingen_Yung_2021Lorentz}).

Somewhat underpinning all these is a certain interplay between $L^p$ and weak-$L^p$ (also known as $L^{p,\infty}$ in the scale of Lorentz spaces). 
In \Cref{sect:Lorentz} we review the basics about Lorentz spaces and some related facts about real interpolation; we also explain a simple instance of how an $L^p$ norm on a lower dimensional space can be realized as a weak-$L^p$ norm in a higher dimensional space, a phenomenon that drives our main theorems. 
In \Cref{sect:BBM} we further motivate our main results by looking at difference quotient characterizations of fractional Sobolev spaces (also known as diagonal Besov spaces) and the related BBM formula. 
In \Cref{sect:diffquot} we state and comment on the proofs of our main results, which is a one-parameter family formulae of Sobolev / $\BV$ norms and $L^p$ norms.
In \Cref{sect:applications} we give some applications towards Gagliardo-Nirenberg interpolation, emphasizing the importance of the existence of a one-parameter family of formulae, and contrasting our formulae for $\BV$ with the one-parameter family of embeddings by Cohen, Dahmen, Daubechies and DeVore.
In \Cref{sect:further} we discuss some other related works and some possible further directions.

The current paper is based on the lecture delivered by the last named author at the award ceremony of the inaugural Antonio Ambrosetti medal.
\Cref{thm6} and some of the material in \Cref{sect:applications} are extensions of known results, which may not have appeared explicitly in the literature before.

\subsubsection*{\it Acknowledgements} 
The idea of replacing strong $L^p$ by weak-$L^p$ in the context of Sobolev-Gagliardo-Nirenberg-type inequalities involving $W^{1,1}$ goes back to the work of L. Greco and R. Schiattarella \cite{MR4106822}.
We are grateful to Carlo Sbordone for calling the attention of one of us (HB) to their result.
P.-L.Y. would like to thank Qingsong Gu for sharing his insights in their collaboration on \cite{MR4249777}.
A.S. and P.-L.Y. would like to thank the  Hausdorff Research Institute of Mathematics and the organizers of the trimester program ``Harmonic Analysis and Analytic Number Theory'' for a pleasant working environment in the summer of 2021. 
The research was supported in part by NSF grant  DMS-2054220 (A.S.) and  by a Future Fellowship FT200100399 from the Australian Research Council (P.-L.Y.).

\section{Lorentz spaces and real interpolation} \label{sect:Lorentz}

Let $(X,\nu)$ be a measure space.
For $1 \leq p < \infty$, if $f \in L^p(\nu)$ for some measure $\nu$, then for every $\lambda > 0$,
\[
\|f\|_{L^p(\nu)}^p = \int_X |f|^p \, \dif\nu \geq \lambda^p \, \nu \{x \colon |f(x)| > \lambda\}.
\]
In particular, if $f \in L^p(\nu)$, then 
\begin{equation}
\label{eq_Oth0beiZau1soo7ne8Chaefo}
\sup_{\lambda > 0} \Big( \lambda \, \nu\{x \colon |f(x)| > \lambda\}^{1/p} \Big) < \infty;
\end{equation}
but the converse is not necessarily true.
If $f$ is measurable on $X$ and the supremum in \eqref{eq_Oth0beiZau1soo7ne8Chaefo} is finite, then $f$ is said to be in weak-$L^p(\nu)$. 
Its weak-$L^p$ (quasi)-norm is defined as the supremum in \eqref{eq_Oth0beiZau1soo7ne8Chaefo}, and denoted by $[f]_{L^{p,\infty}(\nu)}$.
A classic example is given by $f(x) = |x|^{-n/p}$; it is in weak-$L^p(\dif x)$ on $\R^n$, because 
\[
\mathcal{L}^n \{x \in \R^n \colon |x|^{-n/p} > \lambda\} = \mathcal{L}^n \{x \in \R^n \colon |x| \leq \la^{-p/n}\} \simeq \lambda^{-p}.
\]
(Henceforth we write $\mathcal{L}^n$ for Lebesgue measure on $\R^n$.) 
This $f$ is not in $L^p(\dif x)$, because 
\[
\int_{\R^n} |f|^p \, \dif x = \int_{\R^n} |x|^{-n} \, \dif x = +\infty,
\]
and this serves as a motivation of our main result in what follows.

Later on we will also need the Lorentz spaces $L^{p,r}(\nu)$, which for $1 \leq p, r < \infty$ are defined as the space of all measurable $f$ on $X$ with
\[
[f]_{L^{p,r}(\nu)} \coloneqq \Big(r \int_0^{\infty} \lambda^r \nu\{x \colon |f(x)| > \lambda\}^{r/p} \frac{\dif \lambda}{\lambda} \Big)^{1/r} < \infty.
\]
They arise as real interpolation spaces: if $1 \leq p_0 < p_1 \leq \infty$ and $\frac{1}{p} = \frac{1-\theta}{p_0} + \frac{\theta}{p_1}$ for some $0 < \theta < 1$, then for $1 \leq r \leq \infty$,
\[
L^{p,r}(\nu) = [L^{p_0}(\nu),L^{p_1}(\nu)]_{\theta,r}
\]
where for any Banach spaces $B_0$ and $B_1$, and $f \in B_0 + B_1$, the interpolation norm is defined as 
\[
\|f\|_{[B_0,B_1]_{\theta,r}} \coloneqq \Big( \int_0^{\infty} t^{-\theta} \inf_{f = f_0 + f_1} (\|f_0\|_{B_0} + t \|f_1\|_{B_1} )^r \frac{\dif t}{t} \Big)^{1/r}
\]
for $1 \leq r < \infty$, and
\[
\|f\|_{[B_0,B_1]_{\theta,\infty}} \coloneqq \sup_{0 < t < \infty} t^{-\theta} \inf_{f = f_0 + f_1} (\|f_0\|_{B_0} + t \|f_1\|_{B_1} )
\]
when $r = \infty$.
It is also well-known that 
\[
[f]_{L^{p,r}(\nu)} = \|f\|_{L^p(\nu)} \quad \text{if $r = p$}.
\]

In the rest of this article, we will be exploiting a relationship between the $L^p$ norm on a space and the weak-$L^p$ quasi-norm on a higher dimensional space. The simplest instance of this might be the identity
\begin{equation} \label{eq:Tao}
\|f\|_{L^p(X, \, \nu)} = \left[ \frac{f(x)}{y^{1/p}} \right]_{L^{p,\infty}(X \times (0,\infty), \, \nu \, \dif y)}, \quad 1 \leq p < \infty;
\end{equation} 
we thank Terence Tao for communicating this comment. 

\section{The Bourgain-Brezis-Mironescu formula} \label{sect:BBM}

Let's introduce now a shorthand for the first order difference operator: 
\begin{equation} \label{eq:diffop}
\Delta_h u(x) \coloneqq u(x+h)-u(x) \quad \text{for $x, h \in \R^n$}. 
\end{equation}
For $h$ small, and $u \in C^{\infty}$, we have then $|\Delta_h u(x)| \simeq |\nabla u(x) \cdot h|$. Hence roughly speaking, we might believe that $|\nabla u(x)| \simeq \frac{|\Delta_h u(x)|}{|h|}$, at least if one averages over all possible directions for $h$. As a result, to express $\|\nabla u\|_{L^p(\R^n)}$ using a difference quotient instead of a gradient, a naive guess might be to try
\[
\iint_{\R^{2n}} \frac{|\Delta_h u(x)|^p}{|h|^p} \, \dif  h \, \dif  x
\quad \text{
in place of } \quad
\int_{\R^n} |\nabla u(x)|^p \, \dif  x.
\]
This is not working, because if $u_{t}(x) \coloneqq u(t x)$, then 
\[
\int_{\R^n} |\nabla u_{t}(x)|^p \, \dif  x \quad \text{scales like $t^{p-n}$}
\]
but
\[
\iint_{\R^n \times \R^n} \frac{|\Delta_h u_{t}(x)|^p}{|h|^p} \, \dif  h \, \dif  x \quad \text{scales like $t^p$}.
\]
A proper scaling will be achieved if we consider
\[
\iint_{\R^{2n}} \frac{|\Delta_h u(x)|^p}{|h|^p} \frac{\dif  h \, \dif  x}{|h|^n}
\]
instead; if given $b \in \R$ we introduce a difference quotient\footnote{Note that our notation here is different from that in \cite{BSVY}; the $Q_b u(x,y)$ in \cite{BSVY} would be written as $\cQ_{1+b} u(x,h)$ with $h = y-x$ in the current paper. As a result, the set $E_{\la,b}[u]$ in \cite{BSVY} should be compared to the set $\cE_{\la,1+b}[u]$ in \eqref{eq:Esetdef} below.}
\begin{equation} \label{eq:diffquot}
\cQ_b u(x,h) \coloneqq \frac{|\Delta_h u(x)|}{|h|^b},
\end{equation}
then the above suggests that we consider the integral $\iint_{\R^{2n}} \cQ_{1+\frac{n}{p}} u(x,h)^p \, \dif  h \, \dif  x$.

This idea works if we are dealing with fractional Sobolev spaces. 
Indeed, for $0 < s < 1$ and $1 \leq p < \infty$, the fractional Sobolev space $\dot{W}^{s,p}$ is the space of all $u \in L^1_{\loc}(\R^n)$ such that
\[
\begin{split}
\|u\|_{\dot{W}^{s,p}}^p &\coloneqq \iint_{\R^{2n}} \cQ_{s+\frac{n}{p}} u(x,h)^p \, \dif  h \, \dif  x = \iint_{\R^{2n}} \frac{|u(x+h)-u(x)|^p}{|h|^{sp+n}} \, \dif  h \, \dif  x < \infty.
\end{split}
\]
When $1 < p < \infty$, it is known to be equal to the diagonal Besov space $\dot{B}^s_{p,p}$ with comparable norms.
So this suggests again that maybe $\|\nabla u\|_{L^p(\R^n)}$ should be compared to $\|\cQ_{ 1+\frac{n}{p}}u\|_{L^p(\R^{2n}, \, \dif x\, \dif  h)}$? 
Unfortunately this does not work; as observed in \cite{Bourgain_Brezis_Mironescu_2001} (see also \cite{Brezis_2002}), even for $u \in C^{\infty}_c(\R^n)$, unless $u \equiv 0$, the $L^p\!$ norm on $\R^{2n}$ is always infinite! The issue here is that $|h|^{-n/p} \notin L^p(\dif h)$ on $\R^n$. 
On the other hand, we did see $|h|^{-n/p} \in L^{p,\infty}(\dif h)$, and this will partly motivate our main result in the next section.

For now, let's recall the ``BBM-formula'', by Bourgain, Brezis and Mironescu in \cite{Bourgain_Brezis_Mironescu_2001}, which explores what happens to $\|u\|_{\dot{W}^{s,p}}$ as $s \to 1^-$.
On $\R^n$, one instance of this formula says for $1 \leq p < \infty$ and (say) $u \in C^1_c$, that we have\footnote{The original BBM formula was stated and proved for bounded domains in \cite{Bourgain_Brezis_Mironescu_2001}, but the result extends easily to the whole $\R^n$. See e.g.\ \cite{BSY_BBM}*{Appendix A}.}
\begin{equation} \label{eq:BBM}
\lim_{s \to 1^-} (1-s) \|u\|_{\dot{W}^{s,p}}^p = \lim_{s \to 1^-} (1-s) \|\cQ_{ s+\frac{n}{p}}u\|_{L^p(\R^{2n}, \, \dif x\, \dif  h)}^p = \frac{k(p,n)}{p} \|\nabla u\|_{L^p}^p
\end{equation}
where $k(p,n)$ is given explicitly by
\begin{equation} \label{eq:kpn}
k(p,n) \coloneqq  \int_{\S^{n-1}} |e \cdot \omega|^p \, \dif \omega, \quad e \in \S^{n-1}.
\end{equation}
In particular, $\|\cQ_{s+\frac{n}{p}} u\|_{L^p(\R^{2n}, \, \dif x \, \dif  h)}$ blows up like $(1-s)^{-1/p}$ as $s \to 1^-$ unless $u$ is a constant, another indication that $\|\cQ_{1+\frac{n}{p}} u\|_{L^p(\R^{2n}, \, \dif x \, \dif  h)}$ is not good for computing $\|\nabla u\|_{L^p}$.
Our first main result in the next section offers an alternative point of view, that does not involve varying $s$ (as in \eqref{eq:BBM}), but involves a weak-$L^p$ norm instead of the $L^p$ norm on $\R^{2n}$; remember $|h|^{-n/p}$ is not in $L^p(\dif h)$, but it is in weak-$L^p(\dif h)$.
Before we close this section though, we mention a related result of Maz'ya and Shaposhnikova \cite{MR1940355}.
They explored what happens to $\|u\|_{\dot{W}^{s,p}}$ as $s \to 0^+$ instead of $s \to 1^-$, and showed that for (say) $u \in C^1_c$ and $1 \leq p < \infty$,
\begin{equation} \label{eq:MSh}
\lim_{s \to 0^+} s \|u\|_{\dot{W}^{s,p}}^p = \lim_{s \to 0^+} s \|\cQ_{s+\frac{n}{p}}u\|_{L^p(\R^{2n}, \, \dif x \, \dif h)}^p = \frac{2 \sigma_{n-1}}{p} \|u\|_{L^p}^p 
\end{equation}
where $\sigma_{n-1}$ is the surface area of $\S^{n-1}$.
This result will be compared to \Cref{thm6} below.

\section{Difference quotient characterizations} \label{sect:diffquot}

Our first result is a characterization of $\|\nabla u\|_{L^p}$ for $u \in C^{\infty}_c(\R^n)$.
We use the notations of difference and difference quotient introduced in \eqref{eq:diffop} and \eqref{eq:diffquot}.

\begin{thm}[see \cite{Brezis_VanSchaftingen_Yung_2021}] \label{thm1}
Let $n \ge 1$, $1 \leq p < \infty$ and $u \in C^{\infty}_c(\R^n)$. 
Then
\begin{equation} \label{eq:sup_BVY}
\|\nabla u\|_{L^p} \simeq [\cQ_{1+\frac{n}{p}} u]_{L^{p,\infty}(\R^{2n}, \, \dif x \, \dif h)} = \Big[ \frac{\Delta_h u}{|h|^{1+\frac{n}{p}}} \Big]_{L^{p,\infty}(\R^{2n},\, \dif x \, \dif h)}.
\end{equation}
In other words, for $\lambda > 0$ and $b \in \R$, denote by
\begin{equation} \label{eq:Esetdef}
\cE_{\lambda, b}[u] \coloneqq \Big\{(x,h) \in \R^{2n} \colon \cQ_{b} u(x,h) > \lambda  \Big\}
\end{equation}
the superlevel set of $\cQ_{b} u$ at height $\la$. Then
\[
\|\nabla u\|_{L^p}^p \simeq \sup_{\lambda > 0} \Big(\lambda^p \mathcal{L}^{2n}\bigl(\cE_{\lambda, 1+\frac{n}{p}}[u]\bigr) \Big).
\]
In fact, we also have 
\begin{equation} \label{eq:limit_eq_BVY}
\frac{k(p,n)}{n} \|\nabla u\|_{L^p}^p = \lim_{\lambda \to +\infty} \Big(\lambda^p \mathcal{L}^{2n}\bigl(\cE_{\lambda, 1+\frac{n}{p}}[u]\bigr) \Big).
\end{equation}
Here $k(p,n)$ is given by \eqref{eq:kpn}.
\end{thm}

A few remarks are in order.
First, the power $1+\frac{n}{p}$ is dictated by dilation invariance: if $[\cQ_b u]_{L^{p,\infty}(\R^{2n}, \, \dif x \, \dif h)}$ scales like $\|\nabla u\|_{L^p}$ upon replacing $u(x)$ by $u(tx)$ for $t > 0$, then $b = 1+\frac{n}{p}$. 

Next, the limit equality \eqref{eq:limit_eq_BVY} can be proved using Taylor expansion, in a way somewhat reminiscent to the proof of the BBM formula. 
In fact, $\cQ_{1+\frac{n}{p}} u$ is approximately $|h|^{-n/p} |\nabla u(x) \cdot \frac{h}{|h|}|$ when $|h|$ is small, and heuristically $\cE_{\la,1+\frac{n}{p}}[u]$ can be approximated by the set 
\[
\tilde{\cE}_{\la,1+\frac{n}{p}}[u] \coloneqq  \Big\{(x,h) \in \R^{2n} \colon |h|^{-n/p} \Big|\nabla u(x) \cdot \frac{h}{|h|}\Big| > \lambda  \Big\}
\]
when $\la$ is big. 
But for all $\la > 0$,
\[
\lambda^p \mathcal{L}^{2n}\bigl(\tilde{\cE}_{\lambda, 1+\frac{n}{p}}[u]\bigr) = \frac{k(p,n)}{n} \|\nabla u\|_{L^p}^p.
\]
This heuristic can be turned into an actual proof for \eqref{eq:limit_eq_BVY}.

Moreover, in light of the limit equality \eqref{eq:limit_eq_BVY}, to prove \eqref{eq:sup_BVY}, we only need to prove an upper bound for a weak-$L^p$ norm, namely
\[
[\cQ_{1+\frac{n}{p}} u]_{L^{p,\infty}(\R^{2n}, \, \dif x \, \dif h)} \lesssim \|\nabla u\|_{L^p},
\]
which can be done using a Vitali covering lemma (plus the method of rotations in dimensions $n > 1$).
It is somewhat reminiscient of the proof that the Hardy-Littlewood maximal function is bounded from $L^1$ to weak-$L^1$; see also work of Dai, Lin, Yang, Yuan and Zhang \cite{DLYYZ_metricmeasure} who extended our proof to some metric-measure spaces of homogeneous type.

Finally, the weak-$L^p$ quasi-norm in the theorem cannot be replaced by any (bigger) Lorentz $L^{p,r}$ quasi-norm where $r < \infty$.
\begin{thm}[see \cite{Brezis_VanSchaftingen_Yung_2021Lorentz}] \label{thm2}
(i) If $u$ is measurable on $\R^n$ and 
\[
[\cQ_{1+\frac{n}{p}}u]_{L^{p,r}(\R^{2n}, \, \dif x \, \dif h)} < \infty
\]
for some $1 \leq p, r < \infty$, then $u$ is a.e.\ a constant.\\
(ii) Indeed, if $u$ is measurable on $\R^n$ and
\[
\lim_{\lambda \to +\infty} \lambda^p \mathcal{L}^{2n} \Big\{(x,h) \in \R^{2n} \colon \cQ_{1+\frac{n}{p}} u(x,h) > \lambda \Big\} = 0
\]
for some $1 \leq p < \infty$, then $u$ is a.e.\ a constant.
\end{thm}

The difficulty here is that we only know $u$ is measurable. 
If we already know $u \in C^{\infty}_c$, the conclusion of \Cref{thm2} follows from \Cref{thm1}; if we already know $u \in \dot{W}^{1,p}$, then the conclusion follows from our next two theorems.
The case $p = 1$ relies on a result of Poliakovsky \cite{poliakovsky}.

\Cref{thm1} provided a way of computing the $\dot{W}^{1,p}$ norm of a function $u \in C^{\infty}_c(\R^n)$ up to a multiplicative constant.
It turns out there is a natural \emph{one-parameter} family of such formulae for $\|\nabla u\|_{L^p(\R^n)}$, for \emph{general} $u \in \dot{W}^{1,p}$ or $u \in \BV$ (not just for $u \in C^{\infty}_c$); in applications often the question at hand determines which formula one uses within this family.
To state such formulae, let $\ga \in \R$. 
Define the measure 
\begin{equation} \label{eq:nu_def}
\dif \nu_{\ga} \coloneqq |h|^{\ga-n} \, \dif x \, \dif h
\end{equation}
on $\R^{2n}$. 
(The case $\ga = n$ corresponds to the Lebesgue measure $\dif x \, \dif h = \mathcal{L}^{2n}$ we used earlier.)
Then we have two theorems, the first dealing with the case $p > 1$, the second dealing with the case $p = 1$; again we use the constant $k(p,n)$ defined in \eqref{eq:kpn}.

\begin{thm}[see \cite{BSVY}] \label{thm3}
Let $n \ge 1$, $1 < p < \infty$ and $u \in \dot{W}^{1,p}(\R^n)$. 
Then for $\ga \ne 0$, 
\begin{equation} \label{eq:W1p_sup}
\|\nabla u\|_{L^p} \simeq [\cQ_{1+\frac{\ga}{p}} u]_{L^{p,\infty}(\R^{2n}, \, \nu_{\ga})} = \Big[ \frac{\Delta_h u}{|h|^{1+\frac{\ga}{p}}} \Big]_{L^{p,\infty}(\R^{2n}, \, \nu_{\ga})}.
\end{equation}
Furthermore, if $\cE_{\lambda,b}[u]$ is the superlevel set of $\cQ_b u$ at height $\la$ given in \eqref{eq:Esetdef}, then
\begin{equation} \label{eq:W1p_lim}
\frac{k(p,n)}{|\ga|} \|\nabla u\|_{L^p}^p = 
\begin{cases} \lim_{\lambda \to +\infty} \Big(\lambda^p \nu_{\ga}\bigl( \cE_{\lambda,1+\frac{\ga}{p}}[u] \bigr) \Big) \quad &\text{if $\ga > 0$,} \\
\lim_{\lambda \to 0^+} \Big(\lambda^p \nu_{\ga} \bigl( \cE_{\lambda,1+\frac{\ga}{p}}[u] \bigr) \Big) \quad &\text{if $\ga < 0$}.
\end{cases}
\end{equation}
\end{thm}

The case $\ga = -p$ of the above limit equality \eqref{eq:W1p_lim} is due to Nguyen \cite{Nguyen06}.

Next, for the case $p = 1$ we have a similar theorem for $\BV$, but with a number of \emph{additional twists}!

\begin{thm}[see \cite{BSVY}] \label{thm4}
Suppose $n \ge 1$. Then for $\ga \in \R \setminus [-1,0]$ and $u \in \BV(\R^n)$, 
\begin{equation} \label{eq:sup_p=1_good}
\|u\|_{\BV} = \|\nabla u\|_{\cM} \simeq [\cQ_{1+\ga} u]_{L^{1,\infty}(\R^{2n}, \, \nu_{\ga})} =\Big[ \frac{\Delta_h u}{|h|^{1+\ga}} \Big]_{L^{1,\infty}(\R^{2n},\nu_{\ga})}.
\end{equation}
Furthermore, if $\cE_{\lambda,b}[u]$ is the superlevel set of $\cQ_b u$ at height $\la$ given in \eqref{eq:Esetdef}, then the formula
\begin{equation} \label{eq:limit_eq_p=1_good}
\frac{k(1,n)}{|\ga|} \|\nabla u\|_{\cM} = \begin{cases} \lim_{\lambda \to +\infty} \Big(\lambda \nu_{\ga}\bigl(\cE_{\lambda,1+\ga}[u]\bigr) \Big) \quad &\text{if $\ga > 0$} \\
\lim_{\lambda \to 0^+} \Big(\lambda \nu_{\ga}\bigl(\cE_{\lambda,1+\ga}[u]\bigr) \Big) \quad &\text{if $\ga < -1$}
\end{cases}
\end{equation}
holds for $u \in \dot{W}^{1,1}$ but can \emph{fail} for $u \in \BV$ (e.g.\ if $u = \mathbf{1}_{\Omega}$ where $\Omega \subset \R^n$ is any bounded domain with smooth boundary, then the limits above exist but is equal instead to $\frac{k(1,n)}{|\ga+1|} \|\nabla u\|_{\cM}$). 
For $\ga \in [-1,0)$, 
\begin{equation} \label{eq:sup_p=1_bad}
\sup_{u \in C^{\infty}_c(\R^n), \, \|\nabla u\|_{L^1(\R^n)} = 1} [\cQ_{1+\ga} u]_{L^{1,\infty}(\R^{2n}, \, \nu_{\ga})} = +\infty;
\end{equation}
furthermore, the formula 
\begin{equation} \label{eq:limit_eq_p=1_bad}
\frac{k(1,n)}{|\ga|} \|\nabla u\|_{L^1} =
\lim_{\lambda \to 0^+} \Big(\lambda \nu_{\ga}\left(\cE_{\lambda,1+\ga}[u]\right) \Big)
\end{equation}
remains true for all $u \in C^1_c(\R^n)$, but \emph{fails} for $u \in \dot{W}^{1,1}(\R^n)$, and the failure is generic in the sense of Baire category, despite the fact that for all $u \in \dot{W}^{1,1}(\R^n)$ we do have
\begin{equation} \label{eq:upperbdd_by_liminf}
\frac{k(1,n)}{|\ga|} \|\nabla u\|_{L^1} \leq
\liminf_{\lambda \to 0^+} \Big(\lambda \nu_{\ga}\left(\cE_{\lambda,1+\ga}[u]\right) \Big).
\end{equation}
\end{thm}

Note that the range of $\ga$ allowed in \eqref{eq:sup_p=1_good} in \Cref{thm4} is smaller than that in \Cref{thm3}; in fact, \eqref{eq:sup_p=1_bad} shows that \eqref{eq:sup_p=1_good} is false when $\ga \in [-1,0)$, even if one only restricts to functions in $C^{\infty}_c$.
The case $\ga = -1$ of the limiting formula \eqref{eq:limit_eq_p=1_bad} has already been established by Brezis and Nguyen \cite{BN2018}: they have already shown that \eqref{eq:limit_eq_p=1_bad} remains true for all $u \in C^1_c(\R^n)$, but \emph{fails} for general $u \in \dot{W}^{1,1}(\R^n)$.
On the other hand, the counterexamples in the case $\ga \in (-1,0)$ relies on the construction of a Cantor set of dimension $1+\ga$; in fact in this range of $\ga$, one may choose approximations of the associated  Cantor-Lebesgue function (similar to the one in the proof of \Cref{thm7}(ii) below) to establish \eqref{eq:sup_p=1_bad}, and use sums of such functions to construct counterexamples to \eqref{eq:limit_eq_p=1_bad} in $\dot{W}^{1,1}$.

We remark again once $\nu_{\ga}$ is fixed (by fixing $\ga$), the powers $1+\frac{\ga}{p}$ and $1+\ga$ in the denominators of the difference quotients in \Cref{thm3} and \Cref{thm4} are dictated by dilation invariance. 
Also, the two theorems do not address what happens if $\ga = 0$; it turns out things fail strikingly in the case $\ga = 0$.
It can be shown \cite{BSVY}*{Theorem 1.5} that if both $u, \nabla u \in L^1_{\loc}(\R^n)$ then 
\[
\inf \{ \la > 0 \colon \nu_0(\cE_{\la,1}[u]) < \infty\} = \|\nabla u\|_{L^{\infty}};
\]
in particular, if in addition $[\cQ_1 u]_{L^{p,\infty}(\R^{2n},\nu_0)} < \infty$ for some $1 \leq p < \infty$, then $u$ is a.e.\! a constant.

It may be fitting to comment a bit on the proofs of the positive results in \Cref{thm3} and \Cref{thm4}.
Let $n \geq 1$ and $1 \leq p < \infty$.
The issue here is that $C^1_c(\R^n)$ is dense in $\dot{W}^{1,p}(\R^n)$ only when $n > 1$ and $p \geq 1$, or when $n = 1$ and $p > 1$. 
In other words, $C^1_c(\R)$ is not dense in $\dot{W}^{1,1}(\R)$. 
Fortunately it is always possible to approximate a general function in $\dot{W}^{1,p}(\R^n)$ in norm by $C^1$ functions whose \emph{gradients} are compactly supported.
Let $\mathfrak{C}$ denote this latter set of functions.
We have 
\[
C^1_c \subsetneq \mathfrak{C} \subsetneq \dot{W}^{1,p},
\]
and it is no harder to prove, for $u \in \mathfrak{C}$ than for $u \in C^1_c$, the upper bound for $[\cQ_{1+\frac{\ga}{p}} u]_{L^{p,\infty}(\nu_{\ga})}$ for all $\ga \in \R \setminus \{0\}$, or the required upper bound for $\limsup_{\la \to +\infty} \bigl(\lambda^p \nu_{\ga}\bigl( \cE_{\lambda,1+\frac{\ga}{p}}[u] \bigr) \bigr)$ when $\ga > 0$. 
Thus we can pass to limits and conclude the same for a general function in $\dot{W}^{1,p}(\R^n)$ (an additional argument allows us to conclude the same upper bound for $[\cQ_{1+\ga} u]_{L^{1,\infty}(\nu_{\ga})}$ for all $u \in \BV$ and all $\ga \in \R \setminus \{0\}$).
On the other hand, when $\ga < 0$, it is not too hard to establish the desired upper bound for $\limsup_{\la \to 0^+} \bigl(\lambda^p \nu_{\ga}\bigl( \cE_{\lambda,1+\frac{\ga}{p}}[u] \bigr) \bigr)$ for all $u \in C^1_c(\R^n)$. If in addition $p > 1$, or $n \geq 2$ and $\ga < -1$, then we can use the density of $C^1_c(\R^n)$ in $\dot{W}^{1,p}(\R^n)$ mentioned above, together with the upper bound for $[\cQ_{1+\frac{\ga}{p}} u]_{L^{p,\infty}(\nu_{\ga})}$ already proven, to pass to limit to obtain the same conclusion for general $u \in \dot{W}^{1,p}(\R^n)$. The remaining case is then $n = p = 1$ and $\ga < -1$; in this case, one needs to first establish the desired upper bound for $\limsup_{\la \to 0^+} \bigl(\lambda \nu_{\ga}\bigl( \cE_{\lambda,1+\ga}[u] \bigr) \bigr)$ for the bigger class $u \in \mathfrak{C}$ before passing to limit to a general $u \in \dot{W}^{1,1}(\R)$.
Finally, one can establish directly, for all $u \in \dot{W}^{1,p}(\R^n)$, the desired lower bound for $\liminf_{\la \to +\infty} \bigl(\lambda^p \nu_{\ga}\bigl( \cE_{\lambda,1+\frac{\ga}{p}}[u] \bigr) \bigr)$ if $\ga > 0$, and that for $\liminf_{\la \to 0^+} \bigl(\lambda^p \nu_{\ga}\bigl( \cE_{\lambda,1+\frac{\ga}{p}}[u] \bigr) \bigr)$ if $\ga < 0$.
This completes our discussion of all positive results regarding the limiting formulae in \Cref{thm3} and \Cref{thm4}.

Note that if $p > 1$, the lower bound for $[\cQ_{1+\frac{\ga}{p}} u]_{L^{p,\infty}(\nu_{\ga})}$ for $u \in \dot{W}^{1,p}$ in \eqref{eq:W1p_sup} follows from the limiting formulae \eqref{eq:W1p_lim}.
On the other hand, when $p = 1$, $\ga \in \R \setminus [-1,0]$, the limiting formulae \eqref{eq:limit_eq_p=1_good} can fail for $u \in \BV$. 
Thus to prove the lower bound for $[\cQ_{1+\ga} u]_{L^{1,\infty}(\nu_{\ga})}$ in \eqref{eq:sup_p=1_good} for all $u \in \BV$, one must proceed differently. 
The BBM formula comes to our rescue; in fact, the same argument also proves the following theorem, which \emph{characterizes} $\dot{W}^{1,p}$ ($1 < p < \infty$) and $\BV$:

\begin{thm}[see \cite{BSVY}] \label{thm5}
Let $n \ge 1$, $u \in L^1_{\loc}(\R^n)$, $\ga \in \R$. If $[\cQ_{1+\frac{\ga}{p}}u]_{L^{p,\infty}(\R^{2n}, \, \nu_{\ga})} < \infty$, then 
\[
u \in \begin{cases}
\dot{W}^{1,p}(\R^n) \quad & \text{if $1 < p < \infty$} \\
\BV(\R^n) \quad & \text{if $p = 1$}.
\end{cases}
\]
\end{thm}

See also Poliakovsky \cite{poliakovsky}*{Theorem 1.3}, who proved, among other things, the same result for $\ga = n$ under an additional hypothesis $u \in L^p(\R^n)$; in fact, in that case the hypothesis $[\cQ_{1+\frac{n}{p}}u]_{L^{p,\infty}(\R^{2n}, \, \nu_{n})} < \infty$ can be weakened to
\[
\limsup_{\lambda \to +\infty} \Big(\lambda^p \mathcal{L}^{2n} \bigl( \cE_{\lambda,1+\frac{n}{p}}[u] \bigr) \Big) < \infty.
\]

It may be instructive to contrast \Cref{thm5} with \Cref{thm3} and \Cref{thm4}: note that \Cref{thm3} and \Cref{thm4} do not address what happens unless $u \in \dot{W}^{1,p}$ or $u \in \BV$.

To summarize, \Cref{thm3}, \Cref{thm4} and \Cref{thm5} imply that for $u \in L^1_{\loc}(\R^n)$, $1 < p < \infty$ and $\ga \ne 0$, 
\[
u \in \dot{W}^{1,p} \quad \Longleftrightarrow \quad \Big[ \frac{\Delta_h u}{|h|^{1+\frac{\ga}{p}}} \Big]_{L^{p,\infty}(\R^{2n},\nu_{\ga})}
= [\cQ_{1+\frac{\ga}{p}}u]_{L^{p,\infty}(\R^{2n}, \, \nu_{\ga})} 
< \infty.
\] 
Similarly, for $u \in L^1_{\loc}(\R^n)$ and $\ga \in \R \setminus [-1,0]$, 
\[
u \in \BV  \quad \Longleftrightarrow \quad 
\Big[ \frac{\Delta_h u}{|h|^{1+\ga}} \Big]_{L^{1,\infty}(\R^{2n},\nu_{\ga})}
= [\cQ_{1+\ga}u]_{L^{1,\infty}(\R^{2n}, \, \nu_{\ga})} 
< \infty.
\] 

In a slightly different direction, in place of $\|\nabla u\|_{L^p(\R^n)}$, one can also obtain a similar one parameter family of formulae for $\|u\|_{L^p(\R^n)}$. 
\begin{thm} \label{thm6}
Let $n \ge 1$, $1 \leq p < \infty$ and $u \in L^p(\R^n)$. Then for $\ga \ne 0$,
\begin{equation} \label{eq:Lp_sim}
\|u\|_{L^p} \simeq [\cQ_{\frac{\ga}{p}} u]_{L^{p,\infty}(\R^{2n}, \, \nu_{\ga})} = \Big[ \frac{\Delta_h u}{|h|^{\frac{\ga}{p}}} \Big]_{L^{p,\infty}(\R^{2n}, \, \nu_{\ga})}.
\end{equation}
Furthermore, if $\cE_{\lambda,b}[u]$ is the superlevel set of $\cQ_b u$ at height $\la$ given in \eqref{eq:Esetdef}, then
\begin{equation} \label{eq:Lp_limit}
\frac{2\sigma_{n-1}}{|\ga|} \|u\|_{L^p}^p = 
\begin{cases} 
\lim_{\lambda \to 0^+} \bigl(\lambda^p \nu_{\ga}\bigl(\cE_{\lambda,\frac{\ga}{p}}[u]\bigr) \bigr) \quad &\text{if $\ga > 0$} \\
\lim_{\lambda \to +\infty} \bigl(\lambda^p \nu_{\ga}\bigl(\cE_{\lambda,\frac{\ga}{p}}[u]\bigr) \bigr) \quad &\text{if $\ga < 0$}.
\end{cases}
\end{equation}
where $\sigma_{n-1}$ is the surface area of $\mathbb{S}^{n-1}$.
\end{thm}

In this limiting formula \eqref{eq:Lp_limit}, we let $\la \to 0^+$ if $\ga > 0$, and let $\la \to +\infty$ if $\ga < 0$, contrarily to what happened in \Cref{thm3} and \Cref{thm4}.
Also, in the limiting formulae in \Cref{thm3} and \Cref{thm4}, we had a constant $k(p,n)/|\ga|$, and here we had a constant $2\sigma_{n-1}/|\ga|$; these should be compared, respectively, to the constant $k(p,n)/p$ in the BBM formula \eqref{eq:BBM}, and the constant $2\sigma_{n-1}/p$ in the Maz'ya--Shaposhnikova formula \eqref{eq:MSh}. 
The case $\ga = n$ of \eqref{eq:Lp_limit} was proved in \cite{MR4249777}.
Note that we do not obtain a characterization of $L^p(\R^n)$, contrarily to \Cref{thm5}: the $L^{p,\infty}(\nu_{\ga})$ norms of $\cQ_{\frac{\ga}{p}}u$ are finite (in fact zero) when $u$ is a non-zero constant.
We also note that the differences $\Delta_h u(x)$ in \Cref{thm6} can be replaced by other expressions, such that the sums $S_h u(x) \coloneqq u(x+h)+u(x)$, as we will see in the proof below.

\begin{proof}[Proof of \Cref{thm6}]
We consider two cases. 

\textbf{Case 1: Suppose $\ga > 0$.} 
In this case, to prove the upper bound in \eqref{eq:Lp_sim}, note that
\begin{align*}
&\Big\{(x,h) \colon \frac{|\Delta_h u(x)|}{|h|^{\ga/p}} > \la \Big\}\\
\subset & \Big\{(x,h) \colon |h|^{\ga/p} < \frac{2|u(x)|}{\la}\Big\} \bigcup \Big\{(x,h) \colon |h|^{\ga/p} < \frac{2|u(x+h)|}{\la} \Big\}   
\end{align*}
so for any $\la > 0$,
\[
\nu_{\ga} \Big\{(x,h) \colon \frac{|\Delta_h u(x)|}{|h|^{\ga/p}} > \la \Big\}
\leq 2 \nu_{\ga} \Big\{(x,h) \colon |h|^{\ga/p} < \frac{2|u(x)|}{\la} \Big\} 
\lesssim \frac{1}{\la^p} \int_{\R^n} |u(x)|^p \, \dif x.
\]

To prove \eqref{eq:Lp_limit}, and hence the lower bound in \eqref{eq:Lp_sim}, first assume $u$ has compact support in $B_R(0)$. Then
\begin{align*}
\Big\{(x,h) \colon \frac{|\Delta_h u(x)|}{|h|^{\ga/p}} > \la \Big\}
=& \Big\{(x,h) \colon |x| \leq R, |h| > 2R, |h|^{\ga/p} < \frac{|u(x)|}{\la}\Big\} \\
& \quad \bigcup \Big\{(x,h) \colon |x+h| \leq R, |h| > 2R, |h|^{\ga/p} < \frac{|u(x+h)|}{\la}\Big\}  \\
& \quad \bigcup \Big\{(x,h) \colon |h| \leq 2R, |h|^{\ga/p} < \frac{|\Delta_h u(x)|}{\la}\Big\}  
\end{align*}
where all three sets are disjoint. We have
\begin{align*}
& \nu_{\ga} \Big\{(x,h) \colon |x+h| \leq R, |h| > 2R, |h|^{\ga/p} < \frac{|u(x+h)|}{\la}\Big\} \\
= \, & \nu_{\ga} \Big\{(x,h) \colon |x| \leq R, |h| > 2R, |h|^{\ga/p} < \frac{|u(x)|}{\la}\Big\} \\
= \, & \frac{\sigma_{n-1}}{\ga} \int_{|x| \leq R} \Big( \frac{|u(x)|^p}{\la^p} - (2R)^{\ga} \Big)_+ \, \dif x
\end{align*}
so
\begin{align*}
& \lim_{\la \to 0^+} \la^p \nu_{\ga} \Big\{(x,h) \colon |x+h| \leq R, |h| > 2R, |h|^{\ga/p} < \frac{|u(x+h)|}{\la}\Big\} \\
= & \lim_{\la \to 0^+} \la^p \nu_{\ga} \Big\{(x,h) \colon |x| \leq R, |h| > 2R, |h|^{\ga/p} < \frac{|u(x)|}{\la}\Big\} \\
= & \lim_{\la \to 0^+} \frac{\sigma_{n-1}}{\ga} \int_{|x| \leq R} \Big( |u(x)|^p - \la^p (2R)^{\ga} \Big)_+ \, \dif x = \frac{\sigma_{n-1}}{\ga} \int_{\R^n} |u(x)|^p \, \dif x.
\end{align*}
We also have 
\begin{align*}
& \limsup_{\la \to 0^+} \la^p \nu_{\ga} \Big\{(x,h) \colon |h| \leq 2R, |h|^{\ga/p} < \frac{|\Delta_h u(x)|}{\la}\Big\} \\
\leq  & \limsup_{\la \to 0^+} \la^p \nu_{\ga} \Big\{(x,h) \colon |x| \leq R, |h| \leq 2R \Big\} = 0.
\end{align*}
Together this establishes \eqref{eq:Lp_limit} when $u$ has compact support.

If now $u$ is a general $L^p$ function, we approximate by a sequence of compactly supported functions $u_j$, so that $\|u_j - u\|_{L^p(\R^n)} \to 0$ as $j \to +\infty$. Then for any $\varepsilon > 0$, $j \geq 1$,
\[
\Big\{(x,h) \colon \frac{|\Delta_h u(x)|}{|h|^{\ga/p}} > \la \Big\} \subset
 \Big\{(x,h) \colon \frac{|\Delta_h u_j(x)|}{|h|^{\ga/p}} > \la(1-\varepsilon) \Big\} \bigcup \Big\{(x,h) \colon \frac{|\Delta_h (u_j-u)(x)|}{|h|^{\ga/p}} > \la \varepsilon \Big\}
\]
so using the previous result for $u_j$, 
\[
\limsup_{\la \to 0^+} \la^p \nu_{\ga} \Big\{(x,h) \colon \frac{|\Delta_h u(x)|}{|h|^{\ga/p}} > \la \Big\} \leq \frac{1}{(1-\varepsilon)^p} \int_{\R^n} |u_j(x)|^p \, \dif x + \frac{C}{\varepsilon^p} \int_{\R^n} |u_j(x)-u(x)|^p \, \dif x.
\]
Similarly, 
\[
\Big\{(x,h) \colon \frac{|\Delta_h u(x)|}{|h|^{\ga/p}} > \la \Big\} \supset
 \Big\{(x,h) \colon \frac{|\Delta_h u_j(x)|}{|h|^{\ga/p}} > \la(1+\varepsilon) \Big\} \setminus \Big\{(x,h) \colon \frac{|\Delta_h (u_j-u)(x)|}{|h|^{\ga/p}} > \la \varepsilon \Big\}
\]
so
\[
\liminf_{\la \to 0^+} \la^p \nu_{\ga} \Big\{(x,h) \colon \frac{|\Delta_h u(x)|}{|h|^{\ga/p}} > \la \Big\} \geq \frac{1}{(1+\varepsilon)^p} \int_{\R^n} |u_j(x)|^p \, \dif x - \frac{C}{\varepsilon^p} \int_{\R^n} |u_j(x)-u(x)|^p \, \dif x.
\]
We let $j \to +\infty$ before letting $\varepsilon \to 0^+$ in these inequalities to obtain \eqref{eq:Lp_limit}.

\textbf{Case 2: Suppose $\ga < 0$.} 
To prove the upper bound in \eqref{eq:Lp_sim}, note that
\begin{align*}
&\Big\{(x,h) \colon \frac{|\Delta_h u(x)|}{|h|^{\ga/p}} > \la \Big\}\\
\subset & \Big\{(x,h) \colon |h|^{|\ga|/p} > \frac{\la}{2|u(x)|}\Big\} \bigcup \Big\{(x,h) \colon |h|^{|\ga|/p} > \frac{\la}{2|u(x+h)|}\Big\}   
\end{align*}
so for any $\la > 0$,
\[
\nu_{\ga} \Big\{(x,h) \colon \frac{|\Delta_h u(x)|}{|h|^{\ga/p}} > \la \Big\}
\leq 2 \nu_{\ga} \Big\{(x,h) \colon |h|^{|\ga|/p} > \frac{\la}{2|u(x)|}\Big\} 
\lesssim \frac{1}{\la^p} \int_{\R^n} |u(x)|^p \, \dif x.
\]

To prove \eqref{eq:Lp_limit}, and hence the lower bound in \eqref{eq:Lp_sim}, first assume $u \in L^{\infty}$, say $|u| \leq M$, with compact support in $B_R(0)$. Then for $\la > 2M (2R)^{|\ga|/p}$, we have
\[
\frac{|\Delta_h u(x)|}{|h|^{\ga/p}} > \la \quad \Longrightarrow \quad 2M |h|^{|\ga|/p} > 2M (2R)^{|\ga|/p} \quad \Longrightarrow \quad |h| > 2R,
\]
in which case at most one of $x, x+h$ can be in $B_R(0)$. So
\begin{align*}
\Big\{(x,h) \colon \frac{|\Delta_h u(x)|}{|h|^{\ga/p}} > \la \Big\}
=& \Big\{(x,h) \colon |x| \leq R, |h| > 2R, |h|^{|\ga|/p} > \frac{\la}{|u(x)|}\Big\} \\
& \quad \bigcup \Big\{(x,h) \colon |x+h| \leq R, |h| > 2R, |h|^{|\ga|/p} > \frac{\la}{|u(x+h)|}\Big\}   
\end{align*}
where the two sets in the union are disjoint. Hence
\begin{align*}
\nu_{\ga} \Big\{(x,h) \colon \frac{|\Delta_h u(x)|}{|h|^{\ga/p}} > \la \Big\}
&= 2 \int_{|x| \leq R} \int_{|h| > \max\{2R, (\frac{\la}{|u(x)|})^{p/|\ga|}\}} |h|^{\ga-n} \, \dif h \, \dif x\\
&= \frac{2 \sigma_{n-1}}{|\ga|} \int_{|x| \leq R} \min\Big\{(2R)^{-|\ga|}, \frac{|u(x)|^p}{\la^p}\Big\} \, \dif x 
\end{align*}
which says
\[
\begin{split}
\la^p \nu_{\ga} \Big\{(x,h) \colon \frac{|\Delta_h u(x)|}{|h|^{\ga/p}} > \la \Big\} 
&= \frac{2 \sigma_{n-1}}{|\ga|} \int_{|x| \leq R} \min\Big\{\la^p (2R)^{-|\ga|}, |u(x)|^p\Big\} \, \dif x \\
&\to \frac{2 \sigma_{n-1}}{|\ga|} \int_{|x| \leq R} |u(x)|^p \, \dif x
\end{split}
\]
as $\la \to +\infty$ by monotone convergence.

If now $u$ is a general $L^p$ function, we approximate by a sequence of bounded, compactly supported functions $u_j$, so that $\|u_j - u\|_{L^p(\R^n)} \to 0$ as $j \to +\infty$. Then for any $\varepsilon > 0$, $j \geq 1$,
\[
\Big\{(x,h) \colon \frac{|\Delta_h u(x)|}{|h|^{\ga/p}} > \la \Big\} \subset
 \Big\{(x,h) \colon \frac{|\Delta_h u_j(x)|}{|h|^{\ga/p}} > \la(1-\varepsilon) \Big\} \bigcup \Big\{(x,h) \colon \frac{|\Delta_h (u_j-u)(x)|}{|h|^{\ga/p}} > \la \varepsilon \Big\}
\]
so using the previous result for $u_j$, 
\[
\limsup_{\la \to \infty} \la^p \nu_{\ga} \Big\{(x,h) \colon \frac{|\Delta_h u(x)|}{|h|^{\ga/p}} > \la \Big\} \leq \frac{1}{(1-\varepsilon)^p} \int_{\R^n} |u_j(x)|^p \, \dif x + \frac{C}{\varepsilon^p} \int_{\R^n} |u_j(x)-u(x)|^p \, \dif x.
\]
Similarly, 
\[
\Big\{(x,h) \colon \frac{|\Delta_h u(x)|}{|h|^{\ga/p}} > \la \Big\} \supset
 \Big\{(x,h) \colon \frac{|\Delta_h u_j(x)|}{|h|^{\ga/p}} > \la(1+\varepsilon) \Big\} \setminus \Big\{(x,h) \colon \frac{|\Delta_h (u_j-u)(x)|}{|h|^{\ga/p}} > \la \varepsilon \Big\}
\]
so
\[
\liminf_{\la \to \infty} \la^p \nu_{\ga} \Big\{(x,h) \colon \frac{|\Delta_h u(x)|}{|h|^{\ga/p}} > \la \Big\} \geq \frac{1}{(1+\varepsilon)^p} \int_{\R^n} |u_j(x)|^p \, \dif x - \frac{C}{\varepsilon^p} \int_{\R^n} |u_j(x)-u(x)|^p \, \dif x.
\]
We let $j \to +\infty$ before letting $\varepsilon \to 0^+$ in these inequalities to obtain \eqref{eq:Lp_limit}.
\end{proof}

\section{Applications to Gagliardo-Nirenberg interpolation} \label{sect:applications}

The existence of a one-parameter family of characterizations in the previous section is not just natural, but useful in applications.
For instance, Cohen, Dahmen, Daubechies and DeVore \cite{CDDD} proved that for any $0 < t < 1$ and any $1 < q < \infty$, if 
\begin{equation} \label{eq:CDDDassump}
t < \tfrac{1}{q},
\end{equation}
and if
$(\frac{1}{p},s) = (1-\theta) (\frac{1}{q},t) + \theta (1,1)$
for some $0 < \theta < 1$, then for any $u \in \BV \cap \dot{W}^{t,q}$, one has the interpolation inequality
\begin{equation} \label{eq:CDDD}
\|u\|_{\dot{W}^{s,p}} \lesssim \|u\|_{\dot{W}^{t,q}}^{1-\theta} \|u\|_{\BV}^{\theta}.
\end{equation}
\begin{center}
\begin{tikzpicture} [scale = 4]
\draw[arrows=->] (0,0)--(1.1,0);
\draw[arrows=->] (0,0)--(0,1.1);
\node at (0,1.2) {{$s$}};
\node at (1.3,0) {{$1/p$}};
\fill [green!50!black] (1,1) circle[radius=0.02cm];
\node at (1.1,1.1) {{{\color{green!50!black} $\BV$}}};
\draw[black,dotted] (0,0)--(1,1);
\fill [green!50!black] (0.7,0.4) circle[radius=0.02cm];
\node at (0.7,0.25) {{{\color{green!50!black} $\dot{W}^{t,q}$}}};
\fill [blue] (0.85,0.7) circle[radius=0.02cm];
\node [anchor=west] at (0.85,0.7) {{{\color{blue} $\dot{W}^{s,p}$}}};
\draw[blue] (0.7,0.4)--(1,1); 
\node at (1.5,0.5) {{{\color{blue} slope $> 1$}}};
\end{tikzpicture}
\end{center}

Their proof uses bounds for coefficients of wavelet expansions of a general function in $\BV(\R^n)$.
Indeed, let $\psi^0 \coloneqq \varphi$ and $\tilde{\psi}^0 \coloneqq \tilde{\varphi}$ be a pair of one-dimensional compactly supported scaling functions which are in duality:
\[
\int_{\R} \varphi(t) \tilde{\varphi}(t-k) \, \dif t = \delta(k), \quad k \in \Z,
\]
where $\delta$ is the Kronecker delta, and let $\psi^1 \coloneqq \psi$, $\tilde{\psi}^1 \coloneqq \tilde{\psi}$ be their corresponding univariate wavelets. 
Define, for any $e \in E \coloneqq \{0,1\}^n \setminus \{(0,0,\dots,0)\}$, 
\[
\tilde{\psi}^e(x) \coloneqq \tilde{\psi}^{e_1}(x_1) \dots \tilde{\psi}^{e_n}(x_n), \quad x = (x_1, \dots, x_n);
\] 
also define, for any $e \in E$ and any dyadic cube $I = 2^{-j} (k + [0,1]^n)$, 
\[
\tilde{\psi}^e_I(x) \coloneqq 2^{jn} \tilde{\psi}^e(2^j x - k).
\]
For any $\ga \in \R$, one can also define a measure $\tilde{\nu}_{\ga}$ on the product of $E$ with the set of all dyadic cubes $\{I\}$, so that 
\[
\tilde{\nu}_{\ga}(\{(e,I)\}) \coloneqq 2^{-j(\ga+n)}
\]
if $e \in E$ and $I$ has side length $\ell(I) = 2^{-j}$.
A result in \cite{CDDD} says that if $u \in \BV(\R^n)$ and 
\[
u^e_I \coloneqq \int_{\R^n} u(x) \tilde{\psi}^e_I(x) \, \dif x,
\]
then for any $\ga \in \R \setminus [-1,0]$, the sequence $(\frac{u^e_I}{\ell(I)^{1+\ga}})$ indexed by $e$ and $I$ is in weak-$\ell^1$ with respect to $\tilde{\nu}_{\ga}$, with 
\begin{equation} \label{eq:CDDD_weakL1}
\Big[\Big(\frac{u^e_I}{\ell(I)^{1+\ga}}\Big)\Big]_{\ell^{1,\infty}(\tilde{\nu}_{\ga})} \lesssim \|u\|_{\BV} \lesssim \Big\|\Big(\frac{u^e_I}{\ell(I)^{1+\ga}}\Big)\Big\|_{\ell^{1}(\tilde{\nu}_{\ga})}.
\end{equation}
Using \eqref{eq:CDDD_weakL1}, a proof of \eqref{eq:CDDD} was given in \cite{CDDD}; indeed a stronger result was proved there, namely
\begin{equation} \label{eq:CDDD_strong}
[\dot{W}^{t,q},\BV]_{\theta,p} = \dot{W}^{s,p}.
\end{equation}

The inequality \eqref{eq:CDDD_weakL1} bears a superficial resemblance to our difference quotient characterization \eqref{eq:sup_p=1_good} for the $\BV$ norm.
Indeed even the proofs are somewhat similar: both relies on covering lemmas in the range $\ga \in \R \setminus [-n,0]$, and the range $\ga \in [-n,-1)$ for \eqref{eq:CDDD} was dealt with in \cite{CDDD} using the coarea formula, while the same range for \eqref{eq:sup_p=1_good} was dealt with in \cite{BSVY} using the method of rotation.
While we did not manage to use \eqref{eq:sup_p=1_good} to recover a proof of \eqref{eq:CDDD_strong}, the characterization \eqref{eq:sup_p=1_good} does allow us to give a simple  alternative proof of \eqref{eq:CDDD}, which we describe as follows.

\begin{proof}[Proof of \eqref{eq:CDDD}]
Let $\ga_0$ be $-1$ times the slope connecting the points $(1,1)$ and $(\tfrac{1}{q},t)$, i.e.
\begin{equation} \label{eq:ga_0def}
\ga_0 \coloneqq -\frac{1-t}{1-\frac{1}{q}}.
\end{equation}
The assumption \eqref{eq:CDDDassump} shows that $\ga_0 < -1$.
Let $u \in \BV \cap \dot{W}^{t,q}$. Our characterization for the $\BV$ norm (see \Cref{thm4}) shows that 
\begin{equation} \label{eq:formula1}
\|u\|_{\BV} \simeq [\cQ_{1+\ga_0} u]_{L^{1,\infty}(\nu_{\ga_0})}.
\end{equation}
On the other hand, 
\begin{equation} \label{eq:formula2}
\|u\|_{\dot{W}^{t,q}} = \|\cQ_{t+\frac{\ga_0}{q}}u\|_{L^{q}(\nu_{\ga_0})}
\end{equation} 
because from \eqref{eq:nu_def}
\[
\Big( \iint_{\R^{2n}} \frac{|\Delta_h u|^{q}}{|h|^{t q+n}} \, \dif x \, \dif h \Big)^{\frac{1}{q}} 
= \Big( \iint_{\R^{2n}} \frac{|\Delta_h u|^{q}}{|h|^{t q+\ga_0}} \, \dif\nu_{\ga_0} \Big)^{\frac{1}{q}}.
\]
Similarly 
\begin{equation} \label{eq:formula3} 
\|u\|_{\dot{W}^{s,p}} = \|\cQ_{s+\frac{\ga_0}{p}}u\|_{L^{p}(\nu_{\ga_0})}.
\end{equation} 
But since $\frac{1}{p} = (1-\theta) \frac{1}{q} + \theta$, we have, for any measurable function $F$, that  
\begin{equation} \label{eq:interpol_Lorentz}
\|F\|_{L^p(\nu_{\ga_0})} \lesssim \|F\|_{L^{q}(\nu_{\ga_0})}^{1-\theta} [F]_{L^{1,\infty}(\nu_{\ga_0})}^{\theta}; 
\end{equation}
indeed, for any $\la > 0$,
\begin{align*}
    \int |F|^p \, \dif \nu_{\ga_0} 
    &= \int_{|F| \geq \lambda} |F|^p \, \dif \nu_{\ga_0}  + \int_{|F| < \lambda} |F|^p \, \dif \nu_{\ga_0} \\
    &\leq \frac{1}{\lambda^{q-p}} \int |F|^q \, \dif \nu_{\ga_0} + \int_0^{\lambda} s^{p-1} \nu_{\ga_0} \{|F| > s\} \, \dif s \\
    &\leq \frac{1}{\lambda^{q-p}} \|F\|_{L^{q}(\nu_{\ga_0})}^q + \frac{\lambda^{p-1}}{p-1} [F]_{L^{1,\infty}(\nu_{\ga_0})},
\end{align*}
so choosing $\lambda$ for which
\[
\frac{1}{\lambda^{q-p}} \|F\|_{L^{q}(\nu_{\ga_0})}^q = \frac{\lambda^{p-1}}{p-1} [F]_{L^{1,\infty}(\nu_{\ga_0})}
\]
we obtain \eqref{eq:interpol_Lorentz}.
We apply \eqref{eq:interpol_Lorentz} to the function $F \coloneqq \cQ_{s+\frac{\ga_0}{p}} u = \cQ_{t+\frac{\ga_0}{q}} u = \cQ_{1+\ga_0} u$; note that our choice of $\ga_0$ ensures $s + \frac{\ga_0}{p} = t + \frac{\ga_0}{q} = 1 + \ga_0$ (they are all equal to the $y$-intercept of the line joining $(1,1)$ and $(\tfrac{1}{q},t)$). 
Using \eqref{eq:formula1}, \eqref{eq:formula2} and \eqref{eq:formula3}, we obtain \eqref{eq:CDDD}, as desired.

We note that the $\ga_0$ we used above when invoking \Cref{thm4} is dictated by the points $(\tfrac{1}{q},t)$ and $(1,1)$, and this proof does not work if we had chosen other values of $\ga_0$.
\end{proof}

The previous proof made crucial use of the assumption $t < \tfrac{1}{q}$ in \eqref{eq:CDDDassump}, because \eqref{eq:formula1} only holds when $\ga_0 \notin \R \setminus [-1,0]$.
In fact as was shown in \cite{MR3813967}, the inequality \eqref{eq:CDDD} does \emph{not} hold when $t \geq \tfrac{1}{q}$. Nevertheless, a simple adaptation of the above proof of \eqref{eq:CDDD} yields part (i) of the following theorem:

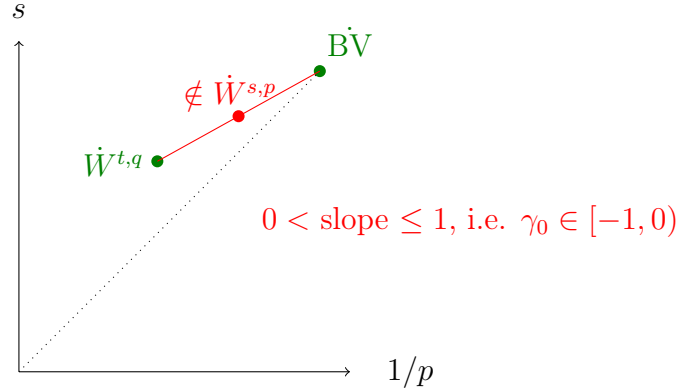
\begin{center}
\begin{tikzpicture} [scale = 4]
\draw[arrows=->] (0,0)--(1.1,0);
\draw[arrows=->] (0,0)--(0,1.1);
\node at (0,1.2) {{$s$}};
\node at (1.3,0) {{$1/p$}};
\fill [green!50!black] (1,1) circle[radius=0.02cm];
\node at (1.1,1.1) {{{\color{green!50!black} $\BV$}}};
\draw[black,dotted] (0,0)--(1,1);
\fill [green!50!black] (0.46,0.7) circle[radius=0.02cm];
\node [anchor=east] at (0.45,0.7) {{{\color{green!50!black} $\dot{W}^{t,q}$}}};
\fill [red] (0.73,0.85) circle[radius=0.02cm];
\node [anchor=south] at (0.7,0.85) {{{\color{red} $\notin \dot{W}^{s,p}$}}};
\draw[red] (0.46,0.7)--(1,1);
\node at (1.5,0.5) {{{\color{red} $0 < $ slope $\leq 1$, i.e. $\ga_0 \in [-1,0)$}}};
\end{tikzpicture}
\end{center}

\begin{thm} \label{thm7}
Let $n \ge 1$, $0 < t < 1$, $1 < q < \infty$. Suppose $t \geq \tfrac{1}{q}$ and $(\frac{1}{p},s) = (1-\theta) (\frac{1}{q},t) + \theta (1,1)$ for some $0 < \theta < 1$. Let $\ga \in \R \setminus [-1,0]$. Then the following hold.
\begin{enumerate}[(i)]
\item Let $r = \tfrac{q}{1-\theta}$. For any $u \in \BV \cap \dot{W}^{t,q}$, one has the interpolation inequality
\begin{equation} \label{eq:BVYinterpol2}
[\cQ_{s+\frac{\ga}{p}}u]_{L^{p,r}(\nu_{\ga})} \lesssim \|u\|_{\dot{W}^{t,q}}^{1-\theta} \|u\|_{\BV}^{\theta}.
\end{equation}
\item The inequality \eqref{eq:BVYinterpol2} fails for some $u \in C^{\infty}_c$ if $r < \frac{q}{1-\theta}$.
\end{enumerate}
\end{thm}

The case $\ga = n$ was already proved in \cite{Brezis_VanSchaftingen_Yung_2021Lorentz}. The proof of the general case is similar, once \Cref{thm4} is established.

\begin{proof}[Proof of \Cref{thm7}]
(i)
Note that since $(\frac{1}{p},s) = (1-\theta) (\frac{1}{q},t) + \theta (1,1)$, we have
\[
s+\frac{\ga}{p} = (1-\theta)(t+\frac{\ga}{q}) + \theta(1+\ga)
\]
for any $\ga \in \R$. In particular,
\[
\cQ_{s+\frac{\ga}{p}}u = (\cQ_{t+\frac{\ga}{q}}u)^{1-\theta} (\cQ_{1+\ga}u)^{\theta}.
\]
It remains to apply H\"{o}lder's inequality for Lorentz spaces: since $(\frac{1}{p},\frac{1}{r}) = (1-\theta) (\frac{1}{q},\frac{1}{q}) + \theta(1,0)$ for $r = \frac{q}{1-\theta}$, we have
\[
[F^{1-\theta} G^{\theta}]_{L^{p,r}(\nu_{\ga})} \lesssim \|F\|_{L^q(\nu_{\ga})}^{1-\theta} [G]_{L^{1,\infty}(\nu_{\ga})}^{\theta}
\]
for any non-negative measurable functions $F$ and $G$.
Applying this to $F = \cQ_{t+\frac{\ga}{q}}u$ and $G = \cQ_{1+\ga}u$, and then invoking \eqref{eq:formula1} and \eqref{eq:formula2} with $\ga \in \R \setminus [-1,0]$ in place of $\ga_0$, yields the desired inequality \eqref{eq:BVYinterpol2}.

\noindent (ii) The optimality of the above choice of $r$ follows the same proof as in \cite{Brezis_VanSchaftingen_Yung_2021Lorentz}*{Lemma 5.1}, which in turn was based on a construction in \cite{MR3813967}; a related example also appeared in \cite{BSVY}*{Proof of Proposition 6.1}. 
We reproduce some of the constructions for the convenience of the readers.

We first consider the case when the dimension $n = 1$. Let $0 < t < 1$ and $1 < q < \infty$. Suppose first $t > \tfrac{1}{q}$. As in \eqref{eq:ga_0def} we define 
$\ga_0 = -\frac{1-t}{1-\frac{1}{q}}$; 
this time $\ga_0 \in (-1,0)$. If $(\tfrac{1}{p},s) = (1-\theta)(\tfrac{1}{q},t) + \theta(1,1)$ for some $0 < \theta < 1$, then $1+\ga_0 = s+\frac{\ga_0}{p}$ whose common value we denote by $\alpha$. Let $\varepsilon := 2^{-1/\alpha} \in (0,1/2)$. Let $g_0$ be an increasing, $C^{\infty}$ function on $\R$ such that $g_0(x) = 0$ for $x < 0$, $g_0(x) = 1$ for $x > 1$. For $j \geq 1$, let $g_j$ be defined on $\R$ by
\[
g_j(x) \coloneqq \frac{1}{2}\left(g_{j-1}(\varepsilon^{-1} x) + g_{j-1}(1-\varepsilon^{-1}(1-x))\right).
\]
The failure of \eqref{eq:BVYinterpol2} when $r < \tfrac{q}{1-\theta}$ can be seen from the inequalities
\begin{equation} \label{eq:g_jBV}
\|g_j\|_{\dot{W}^{1,1}(\R)} = 1,
\end{equation}
\begin{equation} \label{eq:g_jWtq}
\|g_j\|_{\dot{W}^{t,q}(\R)} \lesssim j^{1/q},
\end{equation}
and
\begin{equation} \label{eq:g_jLpr}
[\cQ_{s+\frac{\ga}{p}} g_j]_{L^{p,r}(\R \times \R, \nu_{\ga})} \gtrsim j^{1/r}
\end{equation}
for all $\ga \ne \ga_0$.
To see that these inequalities hold, note that $g_j' \geq 0$, from which \eqref{eq:g_jBV} follows. Also, if 
\[
L_1 = [-(\tfrac{1}{2}-\varepsilon),\tfrac{1}{2}], \quad L_2 = [\tfrac{1}{2},\tfrac{3}{2}-\varepsilon],
\]
then if $|x-y| \geq \tfrac{1}{2}-\varepsilon$ and $g_j(x) \ne g_j(y)$ one must have $(x,y) \in (L_1 \times L_1) \cup (L_2 \times L_2)$ (one can first show if $g_j(x) \ne g_j(y)$ and $|x-y| < \tfrac{1}{2}-\varepsilon$, then both $x$ and $y$ belong to $L_1 \cup L_2$; one can then show that if $(x,y) \in (L_1 \times L_2) \cup (L_2 \times L_1)$ with $|x-y| < \tfrac{1}{2}-\varepsilon$ then $g_j(x) = g_j(y)$). Hence
\begin{align*}
\|g_j\|_{\dot{W}^{t,q}(\R)}^q 
&\leq \|g_j\|_{\dot{W}^{t,q}(L_1)}^q + \|g_j\|_{\dot{W}^{t,q}(L_2)}^q + \iint_{|x-y| \geq \frac{1}{2}-\varepsilon} \frac{1}{|x-y|^{1+t q}} \, \dif x \, \dif y \\
&\leq 2^{1-q} \varepsilon^{1-t q} \|g_{j-1}\|_{\dot{W}^{t,q}(\R)}^q + O(1)
\end{align*}
(we used $t q > 1$ to estimate the last integral) whereas
\[
2^{1-q} \varepsilon^{1-t q} = \varepsilon^{-(1+\ga_0)(1-q)+1-tq} = \varepsilon^{-\ga_0 (\frac{1}{q}-1)q +(1-t)q} = 1;
\]
the estimate \eqref{eq:g_jWtq} now follows by induction on $j$. It remains to establish \eqref{eq:g_jLpr}. We fix $\ga \ne \ga_0$, and define, for $j \ge 0$ and $\la > 0$,
\[
A_{j,\la} \coloneqq \nu_{\ga}\bigl\{(x,h) \in [0,1] \times [0,1] \colon x+h \in [0,1], |\cQ_{s+\frac{\ga}{p}} g_j(x,h)| > \la\bigr\}.
\]
Then for $j \ge 1$ and $\la > 0$, we have
\begin{equation} \label{eq:A_est}
A_{j,\la} \geq \varepsilon^{\ga-\ga_0} A_{j-1,\la \varepsilon^{(\ga-\ga_0)/p} }
\end{equation}
because if $I_1 := [0,\varepsilon]$ and $I_2 := [1-\varepsilon,1]$, then
\begin{align*}
A_{j,\la} & \geq \sum_{i=1}^2 \nu_{\ga}\bigl\{(x,h) \in I_i \times [0,\varepsilon] \colon x+h \in I_i, |\cQ_{s+\frac{\ga}{p}} g_j(x,h)| > \la\bigr\} \\
& = 2 \varepsilon^{1+\ga} \nu_{\ga}\bigl\{(x',h') \in [0,1] \times [0,1] \colon x'+h' \in [0,1], |\cQ_{s+\frac{\ga}{p}} g_{j-1}(x',h')| > 2 \varepsilon^{s+\frac{\ga}{p}} \la\bigr\} \\
& = 2 \varepsilon^{1+\ga} A_{j-1,2 \varepsilon^{s+\frac{\ga}{p}} \la} 
\end{align*}
whereas
\[
2 \varepsilon^{s+\frac{\ga}{p}} = 2 \varepsilon^{s+\frac{\ga_0}{p}} \varepsilon^{\frac{\ga-\ga_0}{p}} = 2 \varepsilon^{\alpha} \varepsilon^{\frac{\ga-\ga_0}{p}} = \varepsilon^{\frac{\ga-\ga_0}{p}},
\]
\[
2 \varepsilon^{1+\ga} = 2 \varepsilon^{1+\ga_0} \varepsilon^{\ga-\ga_0} = 2 \varepsilon^{\alpha} \varepsilon^{\ga-\ga_0} = \varepsilon^{\ga-\ga_0}.
\]
Set $B = B(\ga,\ga_0) := \varepsilon^{-(\ga-\ga_0)}$ so that \eqref{eq:A_est} reads $A_{j,\la} \geq B^{-1} A_{j-1, \la B^{-1/p}}$. Then for $\ell = 1, \dots, j$ and $\la \leq \frac{1}{2}B^{\ell/p}$, we may apply \eqref{eq:A_est} $\ell$~times and invoke $A_{j-\ell,1/2} \gtrsim 1$ to obtain
\begin{equation} \label{eq:g_jintermediate}
\nu_{\ga} \{(x,h) \in [0,1] \times [0,1] \colon x+h \in [0,1], |\cQ_{s+\frac{\ga}{p}} g_j(x,h)| > \la\} \gtrsim B^{-\ell}.
\end{equation}
If $\ga > \ga_0$, then $B > 1$, and hence
\[
[\cQ_{s+\frac{\ga}{p}} g_j]_{L^{p,r}(\nu_{\ga})} \gtrsim \left( \sum_{\ell=1}^j \int_{\frac{1}{2}B^{(\ell-1)/p}}^{\frac{1}{2}B^{\ell/p}} \la^{r-1} B^{-\ell r/p} \dif \la \right)^{1/r}  \simeq j^{1/r}.
\]
If on the other hand $\ga < \ga_0$, then $B < 1$, and hence
\[
[\cQ_{s+\frac{\ga}{p}} g_j]_{L^{p,r}(\nu_{\ga})} \gtrsim \left( \sum_{\ell=1}^j \int_{\frac{1}{2}B^{(\ell+1)/p}}^{\frac{1}{2}B^{\ell/p}} \la^{r-1} B^{-\ell r/p} \dif \la \right)^{1/r} \simeq j^{1/r}.
\]
This proves \eqref{eq:g_jLpr} in either case.

Next, suppose still $n = 1$, and assume $t = \tfrac{1}{q}$ so that $s+\tfrac{\ga}{p} = \frac{1+\ga}{p}$. Then we define instead $g_j(x) := g_0(2^j x) g_0(2^j (2-x))$ where $g_0$ is as above.
We then have $\|g_j\|_{\dot{W}^{1,1}(\R)} = 2$, and $\|g_j\|_{\dot{W}^{t,q}(\R)} \lesssim j^{1/q}$. 
The failure of \eqref{eq:BVYinterpol2} when $r < \tfrac{q}{1-\theta}$ follows once we can show that $[\cQ_{s+\frac{\ga}{p}} g_j]_{L^{p,r}(\R \times \R, \nu_{\ga})} \gtrsim j^{1/r}$ for all sufficiently large $j$ (depending only on $\ga$), which we achieve below by considering the cases $\ga > 0$ and $\ga < -1$ separately.
If $\ga > 0$, then when $1 \leq \la \leq (\tfrac{2^{j-1} \ga}{1+\ga})^{\frac{1+\ga}{p}}$, we have
\begin{align*}
&\nu_{\ga} \bigl\{(x,h) \colon x \leq 0, x+h \geq 2^{-j}, \cQ_{s+\frac{\ga}{p}} g_j(x,h) > \la \bigr\} \\&
\geq \, \int_0^{\la^{-\frac{p}{1+\ga}}} \int_{2^{-j}-h}^{0} h^{\ga-1} \, \dif x \, \dif h 
= \frac{\la^{-p}}{\ga+1}  - \frac{\la^{-p \ga / (1+\ga)}}{\ga 2^j}  \\
&\geq \,  \frac{\la^{-p}}{2(\ga+1)} ,
\end{align*}
the last inequality following from our choice of $\la$. 
It follows that for $j$ sufficiently large, 
\begin{align*}
&[\cQ_{s+\frac{\ga}{p}} g_j]_{L^{p,r}(\nu_{\ga})} 
\\
& \quad\gtrsim \biggl(\int_1^{(\tfrac{2^{j-1} \ga}{1+\ga})^{\frac{1+\ga}{p}}} \la^{r-1} \nu_\gamma\{(x,h) \colon x \leq 0, x+h \geq 2^{-j}, \cQ_{s+\frac{\ga}{p}} g_j(x,h) > \la \}^{r/p} \,  \dif \la \biggr)^{1/r} \\
& \quad\gtrsim \biggl(\int_1^{(\tfrac{2^{j-1} \ga}{ 1+\ga})^{\frac{1+\ga}{p}}} \la^{-1} \, \dif \la \biggr)^{1/r} \simeq j^{1/r}.
\end{align*}
On the other hand, if $\ga < -1$, then when $(\frac{2^{j-2} \ga}{1+\ga})^{\frac{1+\ga}{p}} \leq \la \leq 2^{-\frac{1}{p}}$
\begin{align*}
&\nu_{\ga} \{(x,h) \colon x \leq 0, x+h \geq 2^{-j}, \cQ_{s+\frac{\ga}{p}} g_j(x,h) > \la \} \\
\geq \, & \int_{\la^{-\frac{p}{1+\ga}}}^1 \int_{2^{-j}-h}^0 h^{\ga-1} \, \dif x \, \dif h 
= \frac{\la^{-p}-1}{|\ga+1|}  - \frac{\la^{-p \ga / (1+\ga)}-1}{|\ga| 2^j} \\
\geq \, & \frac{\la^{-p}}{2|\ga+1|}  - \frac{\la^{-p \ga / (1+\ga)}}{|\ga| 2^j} \geq \frac{\la^{-p}}{4|\ga+1|} .
\end{align*} 
(In the penultimate inequality, we used $\la^{-p}-1 \geq \frac{1}{2} \la^{-p}$ which holds since $\la \leq 2^{-\frac{1}{p}}$; in the last inequality, we used $(\frac{2^{j-2} \ga}{1+\ga})^{\frac{1+\ga}{p}} \leq \la$.)
It follows that for $j$ sufficiently large, 
\begin{align*}
&[\cQ_{s+\frac{\ga}{p}} g_j]_{L^{p,r}(\nu_{\ga})} 
\\&\quad \gtrsim \biggl(\int_{(\frac{2^{j-2} \ga}{1+\ga})^{\frac{1+\ga}{p}}}^{2^{-\frac{1}{p}}} \la^{r-1} \nu_\gamma \bigl\{(x,h) \colon x \leq 0, x+h \geq 2^{-j}, \cQ_{s+\frac{\ga}{p}} g_j(x,h) > \la \bigr\}^{r/p} \, \dif \la \biggr)^{1/r} \\
& \quad\gtrsim \biggl(\int_{(\frac{2^{j-2} \ga}{1+\ga})^{\frac{1+\ga}{p}}}^{2^{-\frac{1}{p}}} \la^{-1} \, \dif \la \biggr)^{1/r} \simeq j^{1/r}.
\end{align*}
This completes our proof of the optimality of $r$ in the case where the dimension $n = 1$.

Finally, to pass to higher dimensions $n > 1$, we define
\[
u_j(x) \coloneqq g_j(x_1) \eta_1(x_1) \dots \eta_1(x_n) \in C^{\infty}_c
\]
where $\eta_1 \in C^{\infty}_c((-1,2))$ takes values in $[0,1]$ and is such that $\eta_1 = 1$ on $(-1/2,3/2)$. Then one has
\[
\|u_j\|_{\dot{W}^{1,1}(\R^n)} \lesssim \|g_j\|_{L^1(\R)} + \|g_j\|_{\dot{W}^{1,1}(\R)} \lesssim 1, 
\]
and if we write $\eta'(x_2,\dots,x_n) := \eta_1(x_2) \dots \eta_1(x_n)$, then
\[
\|u_j\|_{\dot{W}^{t,q}(\R^n)} \lesssim \left( \|g_1 \eta_1\|_{L^q(\R)}^q \|\eta'\|_{\dot{W}^{t,q}(\R^{n-1})}^q + \|g_1 \eta_1\|_{\dot{W}^{t,q}(\R)}^q  \|\eta'\|_{L(\R^{n-1})}^q \right)^{1/q} \lesssim j^{1/q}
\]
Furthermore, the argument in \cite{BSVY}*{Section 6.3}, together with our estimates above for $g_j$, shows that for $j$ sufficiently large (depending only on $\ga$),
\[
[\cQ_{s+\frac{\ga}{p}} u_j]_{L^{p,r}(\R^n \times \R^n, \nu_{\ga})} \gtrsim j^{1/r}.
\]
Hence if \eqref{eq:BVYinterpol2} were to hold for all $u \in C^{\infty}_c$, then $r \geq \frac{q}{1-\theta}$. 
\end{proof}

\section{Related works and further directions} \label{sect:further}

The left hand side of the inequality \eqref{eq:BVYinterpol2} involves the quasi-norm $[\cQ_{s+\frac{\ga}{p}}u]_{L^{p,r}(\nu_{\ga})}$, which arises in a number of different contexts in \cite{dssvy}. 
In fact, let $\{L_k\}_{k \in \Z}$ be an appropriate family of Littlewood-Paley projections, and $\mu_{\ga}$ be the measure on $\R^n \times \Z$ given by 
\[
\int_{\R^n \times \Z} F(x,k) \, \dif \mu_{\ga} \coloneqq \sum_{k \in \Z} 2^{-k\ga} \int_{\R^n} F(x,k) \, \dif x
\]
for all $F \in C_c(\R^n \times \Z)$. 
For $0 < s < 1$, $1 < p < \infty$, $1 \leq r \leq \infty$ and $\ga \in \R$, one defines $\dot \cB^s_p(\ga,r)$ to be the space of all tempered distributions $u$ on $\R^n$ modulo polynomials, for which $[2^{k(s+\frac{\ga}{p})} L_k u]_{L^{p,r}(\mu_{\ga})} < \infty$; the set of all measurable functions $u$ on $\R^n$ for which $[\cQ_{s+\frac{\ga}{p}}u]_{L^{p,r}(\nu_{\ga})}$ is finite can then be identified with the space $\dot \cB^s_p(\ga,r)$, which also arises in e.g.\ Krepkogorski\u{\i} \cite{MR1300132} as interpolation spaces between the fractional Sobolev spaces.
Various embedding and non-embedding results for $\dot \cB^s_p(\ga,r)$ were also established in \cite{dssvy}; an application towards nonlinear approximation was also given there.

In \cite{dominguez-milman}, Dom\'{i}nguez and Milman extended some of the above results for $\dot{W}^{1,p}$ in an abstract framework.
They proved for instance that if $(X,\nu)$ is a $\sigma$-finite measure space, $1 \leq p < \infty$ and $\{T_t\}_{t > 0}$ is a family of sublinear operators on $L^p(X)$, then for all $f \in L^p(X)$ satisfying
\[
\|T_t f - f \|_{L^{\infty}(X)} \lesssim_f t^{1/p} \quad \text{for all $t > 0$},
\]
one has
\[
\lim_{\lambda \to \infty} \bigl( \lambda \, (\nu \times \mathcal{L}^1) (\cE_{\lambda})^{1/p} \bigr) = \|f\|_{L^p(X)},
\]
where 
\[
\cE_{\lambda} := \Big \{(x,t) \in X \times (0,\infty) \colon \frac{|T_t f(x)|}{t^{1/p}} > \lambda \Big\}.
\]
They found an impressive list of applications, ranging from formulae for $\|\Delta u\|_{L^p(\R^n)}$ and $\|\partial_{x_1} \partial_{x_2} u\|_{L^p(\R^2)}$, to relations between $\|f\|_{L^p(\R^n)}$ with level set estimates for spherical averages of $f$ for $p > \frac{n}{n-1}$, to ergodic theory, etc.
In \cite{poliakovsky}, Poliakovsky established some of our earlier results on Lipschitz domains on $\R^n$.
The works of Dai, Lin, Yang, Yuan and Zhang \cites{DLYYZ_ballfunctionspaces,DLYYZ_metricmeasure} contain other generalizations of some of the above results to situations where the gradient of a function on $\R^n$ is in a weighted $L^p$ space (which can then be extrapolated), and the case where $\R^n$ is replaced by some suitable metric measure spaces.

In \cite{dominguez-milman2021}, Dom\'{i}nguez and Milman revisited the BBM and the Maz'ya-Shaposhnikova formulae from the points of view of interpolation and extrapolation, putting in context certain results from \cite{BSY_BBM}.

A number of interesting questions remain regarding this circle of ideas.
For instance, for $\ga < 0$, if $u \in L^1_{\loc}$ and the \(\liminf\) on the right hand side of \eqref{eq:upperbdd_by_liminf} is finite, must it be true that $u \in \BV$?  
If $u \in \BV$, must $\|u\|_{\BV}$ be bounded by a multiple (depending only on $n$ and $\ga$) of the liminf on the right hand side of \eqref{eq:upperbdd_by_liminf}?
Some similar questions remain open for $\ga > 0$ and for $\dot{W}^{1,p}$ in place of $\BV$.
For $\ga \in \R \setminus \{0\}$, can one understand the failure of \eqref{eq:limit_eq_p=1_good} for $u \in \BV$ using the concept of $\Gamma$--convergence? 
A more detailed description of these questions, including discussions of partial results, can be found in \cite{BSVY}*{Section 7}.
Also it is conceivable that one may be able to recover the Sobolev inequality, and its Lorentz refinement, namely
\[
\|u\|_{L^{p^*,p}} \lesssim \|\nabla u\|_{L^p}, \quad 1 \leq p < n, \quad \frac{1}{p^*} = \frac{1}{p} - \frac{1}{n}
\]
out of \Cref{thm4}.

\setlength{\parskip}{0pt}

\end{document}